\documentclass[11pt,reqno]{amsart}
\usepackage{amssymb}
\theoremstyle{plain}
\newtheorem{theorem}{Theorem}[section]
\theoremstyle{remark}
\newtheorem{remark}[theorem]{Remark}

\theoremstyle{plain}
\newtheorem{corollary}[theorem]{Corollary}
\newtheorem{lemma}[theorem]{Lemma}
\newtheorem{proposition}[theorem]{Proposition}

\newtheorem{hypothesis}[theorem]{Hypothesis}
\numberwithin{equation}{section}

\begin{document}
\title[EXPONENTIAL ERGODICITY]
{UNIFORM EXPONENTIAL ERGODICITY OF STOCHASTIC 
DISSIPATIVE SYSTEMS}
\author{B. Goldys}
\address{School of Mathematics, The University of New South 
Wales, Sydney 2052, Australia}
\email{B.Goldys@unsw.edu.au}
\author{B. Maslowski}
\address{
Institute of Mathematics, Academy of Sciences\\
\v Zitna 25, 11567 Praha 1, Czech Republic}
\email{maslow@math.cas.cz}
\thanks{This work was partially supported by the Small 
ARC Grant Scheme and GA\v CR grant no. 201/01/1197}
\begin{abstract}
We study ergodic properties of stochastic dissipative 
systems with additive noise. We show that the 
system is uniformly exponentially ergodic provided 
the growth of nonlinearity at infinity is faster than 
linear. The abstract result is applied to the stochastic reaction 
diffusion equation in $\mathbb R^d$ with $d\le 3$. 
\end{abstract}
\maketitle
\section{Introduction}
In this paper we deal with a semilinear stochastic 
equation 
\begin{equation}\left\{\begin{array}{l}
dX=\left(AX+F(X)\right)dt+\sqrt {Q}dW,\\
X(0)=x\in E,\end{array}
\right.\label{1}\end{equation}
in a separable Banach space $\left(E,\left\|\cdot\right\|\right)$ 
continuously embedded into a separable 
Hilbert space $H$ with the inner product $\left\langle\cdot ,\cdot\right
\rangle$ and the 
norm $\left|\cdot\right|$. We assume that $F:E\to E$ is a nonlinear 
mapping, 
$\left(W_t\right)$ is a standard cylindrical Wiener 
process in $H$ defined on a probability space 
$\left(\Omega ,\mathcal F,\left(\mathcal F_t\right),\mathbb P\right
)$ and 
$Q=Q^{*}\in\mathcal L(H)$ is nonnegative. Under 
the assumptions stated below equation (\ref{1}) has a 
unique solution  which defines a Markov $E$-valued 
process with the transition semigroup 
\[P_t\phi (x)=\mathbb E_x\phi\left(X(t)\right),\]
and moreover, it has a unique invariant measure $\mu$. 
In this paper we provide conditions under which the 
convergence to an invariant measure is uniformly 
ergodic in the following sense: There exist positive constants
$C$ and $\gamma$ such that
\begin{equation}\left\|P^*_t\nu-\mu\right\|_{var}\le Ce^{-\gamma t}\left
\|\nu -\mu\right\|_{var}\le 2Ce^{-\gamma t},\label{i1}\end{equation}
for any Borel probability measure $\nu$ on $E$, where $\left\|
\cdot \right\|_{var}$ denotes the norm of total variation of
measures and $P^*_t$ is the adjoint Markov semigroup (in some
papers $P^*_t\nu$ is also denoted by $\nu P_t$).  
This result is known as the uniform 
exponential ergodicity of the Markov 
process associated with the transition semigroup $\left (P_t\right
)$. 
Note that the convergence in (\ref{i1}) is uniform 
with respect to all initial probability measures. This 
property is rather unusual on a non-compact state space. 
For example, if $F=0$ then (\ref{i1}) never holds. 
However, we assume below that the growth of $F$ is faster 
than linear at infinity (since $\epsilon >0$ in Hypothesis 
\ref{A2}), and it turns out that (\ref{i1}) is satisfied, that 
is $\left\|P^*_t\nu  -\mu\right\|_{var}$ is small for large values of $
t$, even if 
$\nu =\delta_a$ (say) with the $\left\|a\right\|$ arbitrarily
large. 

The strong (variational) convergence of $P^*_t\nu $ to the
invariant measure for stochastic evolution equations has been 
investigated in numerous papers (see \cite{M1}-\cite{MSe}, 
the monograph \cite{dz2} and 
the references therein or the survey paper \cite{MSs}). The
geometric ergodicity (which corresponds to the convergence
(\ref{i1}) where the constant $C$ may depend on the initial
measure $\nu$) was studied in \cite{JR} and \cite{Sh}. 
If the diffusion process $X$ is reversible, then as a 
corollary of (\ref{i1}) we obtain 
\begin{equation}\int_E\left|P_t\phi (x)-\left\langle\phi ,\mu\right
\rangle\right|^2\mu (dx)\le e^{-\gamma t}\int_E|\phi (x)|^2\mu (d
x),\label{i2}\end{equation}
where $\left\langle\phi ,\mu\right\rangle =\int\phi d\mu$. Existence of the spectral gap for 
dissipative system (\ref{1}) and for other infinite dimensional 
Markov processes has been recently an object of intense study, 
see for example \cite{aida}, \cite{dap1}, \cite{dap2}, \cite{dadego}, 
\cite{wang}, \cite{wu}, \cite{roc}. 
\par
We will formulate now the main assumptions of the 
paper. 
\begin{hypothesis}\label{A0}
There exists an operator $A_0$ in $H$ such that $A_0$ is an 
infinitesimal generator of a $C_0$-semigroup ${\bf S}=\left({\rm S}
(t)\right)$ on $H$ 
and $A$ is a part of $A_0$ in $E$, that is 
\[\mbox{\rm dom}(A)=\left\{x\in\mbox{\rm dom}\left(A_0\right)\cap 
E:A_0x\in E\right\},\]
and $A=\left.A_0\right|\mbox{\rm dom}(A)$. Moreover, we assume that $
A$ 
generates a compact $C_0$-semigroup in $E$ (which 
we again denote by ${\bf S}$) and 
\begin{equation}\int_0^Tt^{-\alpha}\left\|S(t)Q^{1/2}\right\|_{HS}^
2dt<\infty ,\label{2}\end{equation}
for certain $\alpha ,T>0$, where $\left\|B\right\|_{HS}$ stands for a 
Hilbert-Schmidt norm of an operator $B\in\mathcal L(H)$. 
\end{hypothesis}
It follows from Hypothesis \ref{A0}\ that the stochastic 
convolution integral 
\[Z(t)=\int_0^tS(t-s)\sqrt {Q}dW(s),\quad t\ge 0,\]
is well defined and has an $H$-continuous version. 
\par
 Our next assumption 
concerns regularity of the process $Z$. 
\begin{hypothesis}\label{A1}
There exists $E$-valued, $E$-continuous 
version of the process $Z$, such that 
\begin{equation}\sup_{t\ge 0}\mathbb E\left\|Z(t)\right\|^2<\infty 
.\label{3}\end{equation}
\end{hypothesis}
Our next hypothesis is basically a condition on the 
nonlinear term $F$. By $\langle\cdot ,\cdot\rangle_{E,E^{*}}$
we denote the duality between $E$ and $E^*$ and by $\partial
\| \cdot \|$ the subdifferential of the norm $\| \cdot \|$. 
\begin{hypothesis}\label{A2}
The mapping $F:E\to E$ is Lipschitz continuous on bounded 
sets and for each $x\in\mbox{\rm dom}(A)$ there exists $x^{*}\in\partial\left
\|x\right\|$ such 
that for some $k_1,k_2,k_3>0$ 
\begin{equation}\left\langle Ax,x^{*}\right\rangle_{E,E^{*}}\le 0
,\label{4}\end{equation}
\begin{equation}\left\langle F(x+y),x^{*}\right\rangle_{E,E^{*}}\le 
-k_1\left\|x\right\|^{1+\epsilon}+k_2\left\|y\right\|^s+k_3,\quad 
y\in E.\label{5}\end{equation}
\end{hypothesis}
The solution to equation (\ref{1}) is defined as an 
$E$-continuous adapted process $X$ satisfying the 
integral equation 
\begin{equation}X(t)=S(t)x+\int_0^tS(t-s)F(X(s))ds+Z(t),\quad t\ge 
0.\label{6}\end{equation}
\begin{proposition}\label{t1}
Assume that Hypotheses \ref{A0}, \ref{A1} and \ref{A2} 
hold.  
Then for each $x\in E$ there exists 
a unique solution $X$ to equation (\ref{1}). Moreover, the
equation (1.1) defines an $E$-valued Markov process in the usual
way. 
\end{proposition}
\begin{proof}
The existence and uniqueness of solutions to (\ref{1}) 
follows immediately from Theorem 7.10 in \cite{dz1}. The 
Markov property may be shown as in \cite{dz1}. 
\end{proof}
Let $\left(P_t\right)$ be the Markov semigroup associated to 
equation (\ref{1}), that is 
\begin{equation}P_t\phi (x)=\mathbb E_x\phi\left(X(t)\right),\quad 
x\in E,t\ge 0,\phi\in\mathcal M(E),\label{7}\end{equation}
where $\mathbb E_x$ denotes the expectation corresponding to the
initial condition $X(0)=x$ and $\mathcal M(E)$ denotes the space of bounded 
measurable functions on $E$. Set 
\[P\left(t,x,\Gamma\right)=P_tI_{\Gamma}(x),\quad x\in E,\quad 
\Gamma\in\mathcal 
B(E),\]
where $\mathcal B(E)$ stands for the Borel $\sigma$-algebra on
$E$. Let $\mathcal P$ be the set of Borel probability 
measures on $E$ and 
let $\left(P_t^{*}\right)$ denote the adjoint Markov semigroup acting on 
measures, i.e., 
\[P_t^{*}\nu (\Gamma )=\int_EP(t,x,\Gamma )\nu (dx),\quad t\ge 0,\quad
\Gamma\in\mathcal B(E),\quad\nu\in\mathcal P.\] 
Recall that an invariant measure 
$\mu\in\mathcal P$ is defined as a stationary point of the 
dynamical system $\left(P_t^{*}\right)$, that is $P_t^{*}\mu =\mu$ for all $
t\ge 0$. 
Further, recall that the Markov semigroup $\left(P_t\right)$ is called 
strongly Feller on $E$ if $P_t(\mathcal M)$$\subset C_b(E)$ for each 
$t>0$ (or, alternatively, if the mapping $x\to P(t,x,\Gamma )$ is 
continuous on $E$ for each $t>0$ and $\Gamma\in\mathcal B(E)$), and 
$\left(P_t\right)$ is called topologically irreducible if $P(t,x,
U)>0$ for 
each $t>0$, $x\in E$ and every open set $U\subset E$. Our last 
assumption is 
\begin{hypothesis}\label{A3}
The Markov semigroup $\left(P_t\right)$ associated to the solution of 
of equation (\ref{1}) is strongly Feller and topologically 
irreducible. 
\end{hypothesis}
In Propositions \ref{9} and \ref{10} below, sufficient conditions for 
the strong Feller property and topological irreducibility 
are expressed in terms of coefficients of equation 
(\ref{1}). Basically, they are reformulations of 
known results from \cite{M1}, \cite{M2} and \cite{MSe} 
(see also the monographs \cite{dz1} and \cite{dz2}). It is 
well known that Hypothesis \ref{A3}\ and the existence of 
an invariant measure $\mu\in\mathcal P$ yield $P_t^{*}\nu\to\mu$ as $
t\to\infty$ in the 
total variation norm for every initial measure 
$\nu\in\mathcal P$ (see e.g. \cite{Se1}). 
\section{Uniform exponential ergodicity and some 
auxiliary results}
\begin{proposition}\label{t6}
Assume Hypotheses \ref{A0}, \ref{A1} and \ref{A2}. Then 
there exists $M>0$ such that 
\begin{equation}\sup_{x\in E}\sup_{t\ge 1}\mathbb E_x\left\|X(t)\right
\|\le M.\label{11}\end{equation}
\end{proposition}
\begin{proof}
Let us note first that in virtue of the Fernique theorem 
Hypothesis \ref{A1}\ implies 
\[\sup_{t\ge 0}\mathbb E\left\|Z(t)\right\|^p<\infty ,\]
for all $p>0$, since the process $Z$ is Gaussian in $E$. 
 For $x\in E$ set $Y^x(t)=X(t)-Z(t)$, where $X$ is the 
solution to (\ref{1}) starting from $X(0)=x$ so that 
\begin{equation}Y^x(t)=S(t)x+\int^t_0F\left(Y^x(s)+Z(s)\right)ds,
\quad t\ge 0.\label{13}\end{equation}
We will prove first that for each $x\in E$ and $p\ge\frac 12$, 
\begin{equation}\sup_{t\le T}\mathbb E\left\|Y^x(t)\right\|^{2p}<
\infty .\label{14}\end{equation}
for arbitrary fixed $T>0$.
In the proof of (\ref{14}) we follow similar proofs 
(see Theorem 7.10 of \cite{dz1} or Lemma 2.2 of 
\cite{gm}), so we we omit some details. For $\alpha >0$ we 
define $R(\alpha )=\alpha\left(\alpha I-A\right)^{-1}$, and 
\begin{equation}Y_{\alpha}^x(t)=R(\alpha )S(t)x+\int_0^tR(\alpha 
)S(t-s)F\left(Y^x(s)+Z(s)\right)ds,\quad t\le T.\label{15}\end{equation}
(note that (\ref{4}) implies contractivity of ${\bf S}$, so $R(\alpha )$ is
well defined for each $\alpha >0$).
It is well known that 
\begin{equation}Y_{\alpha}^x\to Y^x,\quad\frac {dY_{\alpha}^x}{dt}
-AY_{\alpha}^x-F\left(Y_{\alpha}^x+Z(t)\right)=\delta_{\alpha}^x\to 
0,\label{16}\end{equation}
uniformly in $t\le T$ as $\alpha\to\infty$ (cf. p. 201 of \cite{dz1}). Also, 
\[\frac {d^{-}}{dt}\left\|Y_{\alpha}^x(t)\right\|^{2p}=2p\left\|Y_{
\alpha}^x(t)\right\|^{2p-1}\frac {d^{-}}{dt}\left\|Y_{\alpha}^x(t
)\right\|\]
\begin{equation}\le 2p\left\|Y_{\alpha}^x(t)\right\|^{2p-1}\left(
k_2\left\|Z(t)\right\|^s+k_3+\left\|\delta_{\alpha}(t)\right\|\right
),\label{17}\end{equation}
by Hypothesis \ref{A2}. Therefore, for $t\le T$, 
\begin{equation}\left\|Y_{\alpha}^x(t)\right\|^{2p}\le\left\|Y_{\alpha}^
x(0)\right\|^{2p}+\int_0^t2p\left\|Y_{\alpha}^x(u)\right\|^{2p-1}\left
(k_2\left\|Z(u)\right\|^s+k_3+\left\|\delta_{\alpha}^x(u)\right\|\right
)du.\label{18}\end{equation}
Taking $p=\frac 12$ and passing with $\alpha$ to infinity we obtain 
\begin{equation}\left\|Y^x(t)\right\|\le\left\|Y^x(0)\right\|+\int_
0^t\left(k_2\left\|Z(u)\right\|^s+k_3\right)du,\label{19}\end{equation}
and (\ref{14}) follows for $p=\frac 12$. By (\ref{16}) and
(\ref{18})  we can see also that for $t\le T$ 
the norm $\left\|Y_{\alpha}^x(t)\right\|$ is bounded uniformly in $
\alpha$. Hence, 
passing with $\alpha$ to infinity in (\ref{18}) we arrive at 
\begin{equation}\left\|Y^x(t)\right\|^{2p}\le\left\|Y^x(0)\right\|^{
2p}+\int_0^t2p\left\|Y^x(u)\right\|^{2p-1}\left(k_2\left\|Z(u)\right
\|^s+k_3\right)du.\label{20}\end{equation}
Now it is easy to prove (\ref{14}) for arbitrary $p>0$ by 
induction (with the induction step $\frac 12$) using \ref{A1}\ and 
the H\"older inequality on the right hand side of (\ref{20}) 
(cf. Lemma 2.2 of \cite{gm}). Using Hypothesis \ref{A2}\ 
we find that 
\begin{equation}\frac {d^{-}}{dt}\left\|Y_{\alpha}^x(t)\right\|\le 
-k_1\left\|Y_{\alpha}^x(t)\right\|^{1+\epsilon}+k_2\left\|Z(t)\right
\|^s+k_3+\left\|\delta_{\alpha}(t)\right\|,\label{21}\end{equation}
and proceeding as above we obtain 
\begin{equation}\left\|Y^x(t)\right\|\le\left\|Y^x(\tau )\right\|-k_1\int_
\tau ^t\left\|Y^x(u)\right\|^{1+\epsilon}du+k_2\int_\tau ^t\left\|Z(u)\right
\|^sdu+k_3(t-\tau ),\label{22}\end{equation}
for $0\le \tau \le t$, which by the Jensen inequality yields 
\begin{equation}\mathbb E\left\|Y^x(t)\right\|\le\mathbb
E\left\|Y^x(\tau) \right\|
-k_1\int_\tau ^t\left(\mathbb E\left\|Y^x(u)\right\|\right)^{1+\epsilon}
du+C(t-\tau ),\quad t\ge \tau \ge 0,\label{23}\end{equation}
for a certain $C>0$. Note that by (\ref{14}) the random 
variables $\left\|Y^x(t)\right\|$, $t\le T$, are uniformly integrable, hence 
the function 
\[\phi (t)=\mathbb E\left\|Y^x(t)\right\|\]
is continuous. A standard comparison theorem yields 
\begin{equation}\phi (t)\le y(t),\quad t\ge 0,\label{24}\end{equation}
where $y$ solve the equation 
\begin{equation}\left\{\begin{array}{ll}
\dot {y}=-k_1y^{1+\epsilon}+C,&t\ge 0\\
y(0)=\left\|x\right\|.\end{array}
\right.\label{25}\end{equation}
By (\ref{24}) and (\ref{25}) it follows that 
\begin{equation}\mathbb E\left\|Y^x(t)\right\|\le\max\left(\left(\frac {
2C}{k_1}\right)^{1+\epsilon},\left(\frac 2{k_1\epsilon}+2\right)^{
1/\epsilon}\right),\quad t\ge 1,x\in E,\label{26}\end{equation}
which together with Hypothesis \ref{A1}\ completes the 
proof of (\ref{11}). 
\end{proof}
\begin{lemma}\label{t60}
Assume Hypotheses \ref{A0}, \ref{A1} and \ref{A2}. Then 
there exist a compact $K\subset E$ and $\kappa >0$ such that 
\begin{equation}\inf_{x\in E}P\left(2,x,K\right)\ge\kappa .
\label{12}\end{equation}
\end{lemma}
\begin{proof}
{\em Step 1.\/} We will show first that the set of probability laws 
\[\mathcal P(r)=\left\{\mathcal L\left(Y^x(1)+Z(1)\right):\left\|
x\right\|\le r\right\},\]
is relatively compact in $E$ for each $r>0$. Indeed, 
since the semigroup ${\bf S}$ is compact in $E$, the set 
\[K_1=\left\{S(1)y:\left\|y\right\|\le r\right\}\]
is relatively compact in $E$. Moreover, the 
operator 
\[L^2(0,1;E)\ni f\to Tf=\int_0^1S(1-u)f(u)du\in E,\]
where the integral is defined in the Bochner sense, 
is compact. Therefore, putting 
\[\tilde B\left(r_1\right)=\left\{f\in L^2(0,1;E):\left\|f\right
\|_{L^2(0,1;E)}\le r_1\right\},\]
we find that $T\left(\tilde B\left(r_1\right)\right)$ is relatively 
compact in $E$. Let 
\[\Omega\left(r_2\right)=\left\{\omega\in C(0,1;E):\sup_{t\le 1}\left
\|Z(t)\right\|\le r_2\right\}.\]
If $\omega\in\Omega\left(r_2\right)$ and $\left\|x\right\|\le r$
then invoking (\ref{19}) we obtain 
\[\left\|Y^x(t)\right\|\le r+\int_0^t\left(k_2r_2^s+k_3\right
)du=r+k_2r_2^{s}+k_3,\quad t\le 1,\]
and since $F$ is bounded on bounded sets of $E$,
\[\sup_{x\in B(r),\omega\in\Omega\left(r_2\right)}\sup_{t\le 1}\left
\|F\left(Y^x(t)+Z(t)\right)\right\|\le\sup_{y\in B(R)}\left\|F(y)\right
\|<\infty ,\]
where $R=r+r_2+k_2r_2^{s}+k_3$. Let 
\[f_{x,\omega}(t)=F\left(Y^x(t)+Z(t)\right).\]
Then 
\[\mathcal U(r)=\left\{f_{x,\omega}:x\in B(r),\omega\in\Omega\left
(r_2\right)\right\}\subset\tilde B(R)\]
and therefore the set 
$K_2=T\mathcal U(r)$ is relatively compact in $E$. For a 
given $\eta\in (0,1)$ we choose $r_2$ in such a way that 
\[\mathbb P\left(\Omega\left(r_2\right)\right)\ge 1-\frac 12\eta 
.\]
Let $K_3\subset E$ be such a compact set that 
\[\mathbb P\left(Z(1)\in K_3\right)\ge 1-\frac 12\eta ,\]
and let $\Omega_1=\left\{\omega :Z(1)\in K_3\right\}$. 
Finally, let $K(r)=K_1+K_2+K_3$. Then, for $x\in B(r)$, 
\[\mathbb P\left(Y^x_{}(1)+Z(1)\in K_1+K_2+K_3\right)\ge\mathbb P\left
(\Omega\left(r_2\right)\cap\Omega_1\right)\ge 1-\eta .\]
\par\noindent
{\em Step 2 Conclusion.\/} It follows from Step 1 that for each 
$\eta\in (0,1)$ and 
$r>0$ there exists a 
compact set 
$K(r)\subset E$ such that  
\[\inf_{\left\|y\right\|\le r}P\left(1,y,K(r)\right)>1-\eta .\]
Moreover, (\ref{11}) yields the existence of $R>0$ 
such that 
\[P\left(1,x,B\left(R\right)\right)\ge 1-\eta ,\quad x\in E.\]
Then by the Chapman-Kolmogorov equality
\[P\left(2,x,K(r)\right)\ge\int_{B\left(R\right)}P\left(1,y,K(R)\right
)P\left(1,x,dy\right)\ge (1-\eta )^2,\]
which completes the proof of the lemma. 
\end{proof}
Let us recall some basic concepts of Ergodic Theory of 
Markov chains. Let $\left(X_i\right)$ be an $E$-valued Markov chain 
with the transition kernel $P^m(x,\Gamma )$, $m\in\mathbb N$, $x\in 
E$, 
$\Gamma\in\mathcal B(E)$, and let $\phi\ge 0$ be a nontrivial measure 
on $\mathcal B(E)$. The chain $\left(X_i\right)$ is called 
$\phi$-irreducible if 
for each $\Gamma\in\mathcal B(E)$ with $\phi (\Gamma )>0$ we have 
\begin{equation}\sum_{i=1}^{\infty}P^i(x,\Gamma )>0,\quad x\in E.\label{36}\end{equation}
Recall that a set $\Pi\in\mathcal B(E)$ is called a small set if 
there exist $m\in\mathbb N$ and a nontrivial measure $\lambda\ge 
0$ 
such that 
\begin{equation}\inf_{x\in\Pi}P^m(x,\cdot )\ge\lambda (\cdot ).\label{37}\end{equation}
We will need the following result which is an immediate 
consequence of Lemma 2 in \cite{JJ}, see also Theorem 
5.2.2 in \cite{meyn}. 
\begin{lemma}\label{t7}
Let $\left(X_i\right)_{i\in\mathbb N}$ be $\phi$-irreducible. 
Then there exists a small set $\Pi\in\mathcal B\mathbb(E)$ 
such that $\phi (\Pi$$)>0$. 
\end{lemma}
\begin{theorem}\label{t2}
Assume Hypotheses \ref{A0}-\ref{A2} and \ref{A3}. Then there exists an 
invariant measure $\mu\in\mathcal P$ such that for certain 
constants $C>0$, and $\gamma >0$ we have 
\begin{equation}\left\|P_t^{*}\nu -\mu\right\|_{var}\le Ce^{-\gamma 
t}\left\|\nu -\mu\right\|_{var}\le 2Ce^{-\gamma t}\label{9}\end{equation}
for all $t>0$ and $\nu\in\mathcal P$, where $\left\|\cdot\right\|_{
var}$ stands for 
the norm of total variation of measures. 
\end{theorem}
\begin{proof}\ Consider the skeleton chain 
$\left(X_n\right)$, where $X_i=X(i)$ and for a fixed $x_0\in E$ set 
$\phi (\cdot )=P\left(1,x_0,\cdot\right)$. It is well known that by Hypothesis 
\ref{A3}\ the measures $\left\{P(t,x,\cdot ):t>0,x\in E\right\}$ are equivalent 
 hence the chain $\left(X_n\right)$ is $\phi$-irreducible 
and by Lemma \ref{t7}\ there exists a set $\Pi\in\mathcal B(E)$ 
such that 
\begin{equation}P\left(1,x_0,\Pi\right)>0,\label{38}\end{equation}
and 
\begin{equation}\inf_{x\in\Pi}P\left(m,x,\Gamma\right)\ge\lambda 
(\Gamma ),\quad\Gamma\in\mathcal B(E),\label{39}\end{equation}
for some $m\in\mathbb N$ and a nontrivial measure $\lambda$. By 
(\ref{12}) we have 
\[\inf_{x\in E}P\left(m+3,x,\Gamma\right)\ge\inf_{x\in E}\int_{\Pi}
P\left(3,x,dy\right)P\left(m,y,\Gamma\right)\ge\lambda (\Gamma )\inf_{
x\in E}P\left(3,x,\Gamma\right)\]
\[=\lambda (\Gamma )\inf_{x\in E}\int_EP\left(1,y,\Pi\right)P\left
(2,x,dy\right)\ge\lambda (\Gamma )\inf_{x\in E}\int_KP\left(1,y,\Pi\right
)P\left(2,x,dy\right)\]
\begin{equation}\ge\kappa\lambda (\Gamma )\inf_{y\in K}P\left(1,y
,\Pi\right).\label{40}\end{equation}
By (\ref{38}) and the equivalence of transition measures 
$P\left(1,y,\Pi\right)>0$ for all $y\in E$. Since the function $y
\to P\left(1,y,\Pi\right)$ 
is continuous by the Strong Feller Property and $K$ is 
compact, we obtain 
\begin{equation}\inf_{x\in E}P\left(m+3,x,\Gamma\right)\ge\kappa\delta
\lambda (\Gamma ),\quad\Gamma\in\mathcal(E),\label{41}\end{equation}
for a certain $\delta >0$. For $T=m+3$  it follows that 
\begin{equation}P_{t+T}^{*}\nu (\Gamma )=\int_E\int_EP\left(T,y,\Gamma\right
)P\left(t,x,dy\right)\nu (dx)\ge\tilde{\mu }(\Gamma ),\quad t\ge 
0,\nu\in\mathcal P,\label{42}\end{equation}
where $\tilde{\mu }(\cdot )=\kappa\delta\lambda (\cdot )$. Hence, $
\tilde{\mu}$ is a nontrivial lower bound 
measure and it follows that there exists an invariant 
measure $\mu\in\mathcal P$ (see e.g. \cite{La}). To 
prove the exponential convergence, take arbitrary $\delta_1,\delta_2\in\mathcal 
P$ 
and set $\delta =\delta_1-\delta_2$. we will denote by $\zeta^{+}$ and $\zeta^{
-}$ the 
positive and negative part respectively of a signed 
measure $\zeta$. Obviously, we have 
\begin{equation}\eta :=\delta^{+}(E)=\delta^{-}(E)=\frac 12\left\|\delta\right
\|_{var},\label{43}\end{equation}
and without loss of generality we can assume $\eta >0$. 
Then 
\begin{equation}\left\|P_t^{*}\delta\right\|_{var}=\eta\left\|\left(
P_t^{*}\left(\frac 1{\eta}\delta^{+}\right)-\tilde\mu\right)-\left(P_
t^{*}\left(\frac 1{\eta}\delta^{-}\right)-\tilde\mu\right)\right\|_{
var},\quad t\ge 0.\label{44}\end{equation}
Furthermore, by (\ref{42}) the measures $P_T^{*}\left(\frac 1{\eta}
\delta^{+}\right)-\tilde{\mu}$ and 
$P_T^{*}\left(\frac 1{\eta}\delta^{-}\right)-\tilde{\mu}$ are nonnegative, thus 
\begin{equation}\left\|P_T^{*}\left(\frac 1{\eta}\delta^{+}\right)-\tilde
\mu\right\|_{var}=P_T^{*}\left(\frac 1{\eta}\delta^{+}\right)(E)-\tilde{
\mu }(E)=1-\tilde{\mu }(E),\label{45}\end{equation}
and similarly 
\begin{equation}\left\|P_T^{*}\left(\frac 1{\eta}\delta^{-}\right)-\tilde
\mu\right\|_{var}=1-\tilde{\mu }(E),\label{46}\end{equation}
which by (\ref{43}) and (\ref{44}) yields 
\begin{equation}\left\|P_T^{*}\delta\right\|_{var}\le\eta\left(2-2\tilde
\mu (E)\right)=\left(1-\tilde\mu (E)\right)\left\|\delta\right\|_{va
r}.\label{47}\end{equation}
For $q=1-\tilde{\mu }(E)\in (0,1)$ the semigroup property of $P_t^{
*}$ 
yields 
\begin{equation}\left\|P_{nT}^{*}\delta\right\|_{var}\le q^n\left\|\delta\right
\|_{var},\quad n\ge 1,\label{48}\end{equation}
which is the geometric ergodicity for the chain $\left(X_{nT}\right
)$. 
Set $\beta =-\log q$$>0$ and let $[\alpha ]$ stand for the integer part of 
the real number $\alpha$. Then 
\[\left\|P_t^{*}\delta\right\|_{var}\le\left\|P_{t-\left[\frac tT\right
]T}^{*}P_{\left[\frac tT\right]T}^{*}\delta\right\|_{var}\le\left\|P_{\left
[\frac tT\right]T}^{*}\delta\right\|_{var}\]
\begin{equation}\le e^{-\beta\left[\frac tT\right]}\left\|\delta\right
\|_{var}\le e^{-\frac {\beta}Tt}e^{\frac {\beta}T\left(t-\left[\frac 
tT\right]T\right)}\left\|\delta\right\|_{var}\le e^{\beta}e^{-\frac {
\beta}Tt}\left\|\delta\right\|_{var},\label{49}\end{equation}
and (\ref{9}) follows with $\gamma =\frac {\beta}T$ and $C=e^{\beta}$.
\end{proof}
As a corollary of Theorem \ref{t2} we obtain the 
following property. 
\begin{corollary}\label{t4}
Assume Hypotheses \ref{A0}-\ref{A2} and \ref{A3} and let the Markov 
semigroup $\left(P_t\right)$ be symmetric in $L^2(E,\mu )$. Then there 
exist constants $C,\gamma >0$ such that 
\begin{equation}\left\|P_t\phi - <\phi ,\mu >\right\|_{L^2(E,\mu 
)}\le e^{-\gamma t}\left\|\phi -<\phi ,\mu >\right\|_{L^2(E,\mu
)} \le e^{-\gamma t}\left\|\phi \right\|_{L^2(E,\mu )}
\label{10}\end{equation}
for all $\phi\in L^2(E,\mu )$ and $t\ge 0$.
\end{corollary}
\begin{proof}\ The proof follows easily from 
Theorem \ref{t2} and \cite{chen}, Theorem 1.2 (see also 
\cite{rr}).
\end{proof}
\begin{remark}\label{t5}
(i) In \cite{dadego}, Section 6.2, estimate (\ref{10}) is 
obtained essentially for a 
strongly dissipative symmetric system provided $Q$ is 
boundedly invertible. Then, in Section 6.3, an analogue 
of (\ref{10}) is obtained for systems 
with the nonlinearity $F=F_0+F_1$ with $F_0$ strongly 
dissipative and $F_1$ bounded and $Q$ still boundedly 
invertible. However, in the latter case estimate (\ref{10}) 
holds in $L^2\left(E,\mu_0\right)$, where $\mu_0$ is the unique invariant 
measure of equation (\ref{1}) with $F_1=0$. Hence, our 
result and the result from \cite{dadego} are not exactly 
comparable. Finally, let us note that in our case (\ref{10}) 
holds even if $Q$ is not boundedly invertible, provided 
$F=QDG$, where $DG$ is the gradient of the mapping 
$G:E\to\mathbb R$. 
\par
(ii) Due to \cite{chen} the first inequality in (\ref{10}) 
is equivalent (for symmetric 
$\left(P_t\right)$ ) to (\ref{9}) which however must be satisfied only 
if $\nu\ll\mu$, $\frac {d\nu}{d\mu}\in L^2(E,\mu )$, and with
$C$ possibly dependent on 
$\nu$. Thus the statement of Theorem \ref{t2} is essentially 
stronger than the convergence (\ref{10}).
\end{remark}
For the reader's convenience we will amend this section 
with three propositions which are minor modifications of earlier results 
\cite{MSe}, \cite{M2} and \cite{pertex}, in which
Hypothesis \ref{A3} (strong Feller property and irreducibility)
is verified. 
\begin{proposition}\label{t8}
Assume Hypotheses \ref{A0}, \ref{A1}\  and \ref{A2}. Let 
\[Q_t=\int_0^tS(s)QS^{*}(s)ds,\]
and let 
\begin{equation}S(t)(E)\subset Q_t^{1/2}(H),\quad t>0.\label{50}\end{equation}
If there exists a mapping $u\in C(E,H)$ which is bounded on 
bounded sets and such that $F=Q^{1/2}u$ then the solution 
to (\ref{1}) is strong Feller and irreducible i.e. Hypothesis 
\ref{A3}\ holds. In particular, if $Q$ is boundedly invertible  
and $Q^{1/2}\in\mathcal L(E,H)$ then the above conditions hold 
with $u=Q^{-1/2}F$. 
\end{proposition}
\begin{proof}The Strong Feller Property follows from Theorem 3.1 of 
\cite{MSe}, where applicability of the Girsanov Theorem 
to equation (\ref{1}) is also proved. Since (\ref{50}) 
implies topological irreducibility for the linear equation 
($F=0$), the solution to (\ref{1}) is irreducible as well. 
\end{proof}
%\begin{remark}
%Let us note that without loss of generality we may 
%assume that 
%\[\sup_{t>0}\mbox{\rm tr}\left(Q_t\right)<\infty ,\]
%and in that case the operator $Q_{\infty}=\lim_{t\to\infty}Q_t$ is well 
%defined and of trace class. Let $H_{\infty}=\mbox{\rm im}\left(Q_{
%\infty}^{1/2}\right)$. It 
%follows from Hypothesis \ref{A1}\ that $H_{\infty}$ is 
%ontinuously imbedded into $E$. If moreover (\ref{50}) holds 
%then $S(t)(H)\subset\mbox{\rm im}\left(Q_t^{1/2}\right)\subset H_{\infty}$ as well and the semigroup $\mathbf 
%S$ is a 
%compact $C_0$-semigroup 
%when restricted to $H_0$ endowed 
%$with the norm $\left|x\right|_{\infty}=\left|Q_{\infty}^{-1/2}x\right
%|$, see \cite{fock}. It 
%follows easily that $\mathbf S$ is compact in $E$ as well, 
%hence 
%the compactness assumption in Hypothesis \ref{A0}\ may 
%be dropped. Indeed, for any $t>0$ we have 
%$S(t)=S(t-\epsilon )S(\epsilon )$ for a certain $\epsilon >0$ and therefore 
%we find that $S(\epsilon ):E\to H_{\infty}$ is bounded. Finally, again by 
%\cite{fock}, $S(t-\epsilon )$ is compact in $H_{\infty}$. 
%\end{remark}
Proposition \ref{t8} is applicable basically (though not 
exclusively, cf. \cite{MSe}) to the cases when $Q$ is 
boundedly invertible. In the following two statements 
$Q^{-1}$ may be unbounded. 
\begin{proposition}\label{t9}
Let $Q>0$, assume Hypotheses \ref{A0}, \ref{A1}\  and \ref{A2}\ and let 
one of the following conditions be satisfied: either  
\par\noindent
(i) $S(t)H\subset E$ for $t>0$ and $\left\|S(t)\right\|_{H\to E}\le 
q(t)$ with a 
certain $q\in L^2(0,T)$, or 
\par\noindent
(ii) $Q^{1/2}\in\mathcal L(H,E)$ and $\overline {Q^{1/2}(H)}=E$. 
\par\noindent
Then the solution to equation (\ref{1}) is topologically 
irreducible.
\end{proposition}
\begin{proof}
See Propositions 2.7, 2.8 and 2.11 and Lemma 2.6 of \cite{M3}. 
\end{proof}
\begin{proposition}\label{t10}
Assume Hypotheses \ref{A0}, \ref{A1}\  and \ref{A2}. 
Moreover, assume that for each $n\in  {\mathbb N}$ there exists
a $k_n < \infty$ such that
\begin{equation}\left|F(x)-F(y)\right|\le k_n\left|x-y\right|,\quad\left
\|x\right\|+\left\|y\right\|\le n,\label{52}\end{equation}
(that is $F$ is Lipschitz continuous on bounded sets of $E$ 
with respect to the norm in $H$), $S(t)(H)\subset Q_t^{1/2}(H)$ for 
$t>0$, and 
\begin{equation}\int_0^T\left\|Q_t^{-1/2}S(t)\right\|_{\mathcal L
(H)}dt<\infty ,\label{53}\end{equation}
for a certain $T>0$. Then the solution to (\ref{1}) is 
strongly Feller. 
\end{proposition}
\begin{proof}The proof is a simple combination of 
arguments from \cite{M2} and \cite{pertex} so it is only 
sketched. Let $c>0$ be the norm of the embedding $j:E\to H$. 
For $m\ge 1$, and $x\in E$ set 
\begin{equation}F_m(x)=\left\{\begin{array}{lll}
F(x)&\mbox{\rm if}&|x|\le cm,\\
F\left(\frac {cmx}{|x|}\right)&\mbox{\rm if}&|x|>cm.\end{array}
\right.\label{54}\end{equation}
By (\ref{52}) $F_m$ is uniquely extendible to a bounded, 
globally Lipschitz function on $H$ for each $m\in\mathbb N$. 
Therefore, the solution to the equation 
\begin{equation}\left\{\begin{array}{l}
dX_m(t)=\left(AX_m(t)+F_m\left(X_m(t)\right)\right)dt+\sqrt {Q}dW
(t),\\
X_m(0)=x,\end{array}
\right.\label{55}\end{equation}
is strongly Feller for each $m\in\mathbb N$ by 
\cite{pertex}. Since the paths of solutions to (\ref{55}) 
and (\ref{1}) coincide with high probability if $m$ is large, it 
follows easily from the proof of 
Proposition \ref{t6}\ (i) that 
\begin{equation}\lim_{m\to\infty}\sup_{\left\|x\right\|\le R}\left
\|P_m(t,x,\cdot )-P(t,x,\cdot )\right\|_{var}=0,\label{56}\end{equation}
for all $t>0$, $R>0$ where $P_m$ denotes the transition 
kernel associated with (\ref{55}). Hence the solution to 
(\ref{1}) is strongly Feller as well. 
\end{proof}
\begin{remark}\label{t11}
(i) The existence of an invariant measure for (\ref{1}) has 
been proved independently in \cite{dgz} (cf. also 
\cite{gg}) by a method based on a version of the 
Krylov-Bogolyubov argument. 
\par\noindent
(ii) Note that in the proof of Theorem \ref{t2}\ we have 
proved that the whole space $E$ is a small set for the 
chain $\left(X_n\right)$ (cf. (\ref{40}), (\ref{41})). The 
exponential ergodicity of the chain $\left(X_n\right)$ follows also by 
this fact and Theorem 16.2.2 in \cite{meyn}. However, in the
respective part of the proof of Theorem \ref{t2}\ we prove the
exponential ergodicity directly using a simple argument.

(iii) The method used in the paper can be also easily applied 
to some cases of stochastic evolution equations with 
non-additive noise term; basically to the case when the 
diffusion coefficient is bounded and has bounded inverse. 
and the semigroup ${\bf S}$ is exponentially stable. For example, 
if the nonlinear drift term $F$ obeys Hypothesis 
\ref{A2} and the conditions (C1)-(C5) from the paper 
\cite{MSe} are  satisfied the proof of Theorem \ref{t2} 
can be repeated (with obvious modifications).
\end{remark}
\section{Example}
Consider a stochastic parabolic equation
\begin{equation}\left\{\begin{array}{l}
\frac {\partial}{\partial t}u(t,\xi )=\Delta u(t,\xi )+f(t,\xi )+
\eta (t,\xi ),\quad (t,\xi )\in {\mathbb R}_{+}\times D,\\
u(0,\xi )=x(\xi ),\quad\xi\in D,\\
u(t,\xi )=0,\quad (t,\xi )\in {\mathbb R}_{+}\times\partial D,\end{array}
\right.\label{e1}\end{equation}
on a bounded domain $D\subset {\mathbb R}^d,\quad d\le 3,$ with a smooth boundary $
\partial D$, 
where $f:{\mathbb R}\to {\mathbb R}$ is locally Lipschitz and $\eta$ 
symbolically denotes a noise white in time and, in 
general, dependent on the space variable $\xi .$ The system 
(\ref{e1}) is rewritten in a usual manner as an equation of 
the form (1.1) where we put $H=L^2(D),$ $E=C_0(D),$ 
$A_0=\Delta$ with $\mbox{\rm dom}\left(A_0\right)=H_{0^{}}^1(D)\cap 
H^2(D),$ and $F:E\to E$ is defined as 
the superposition operator, $F(y)(\xi ):=f(y(\xi )),$ where 
$y\in E,$ $\xi\in D.$ 
The noise $\eta$ is modelled in the equation (1.1) 
by the Wiener process $W_t$ and the covariance operator 
$Q\in\mathcal L(H),$ formally we have $\eta =^{}Q^{1/2}\frac {dW}{
dt}.$ 
It is well known that the operator $A_0$ generates a 
strongly continuous semigroup on $H$, its part $A$ on 
the space $E$ generates a strongly continuous semigroup 
on $E$ as required in Hypothesis 1.1.
Assume that there exist positive constants $c_1$, $c_2$, $c_3$, $
s$ and 
$\epsilon$ such that
\begin{equation}f(\alpha +\beta )\operatorname{sgn}\alpha\le -c_1|\alpha |^{1+\epsilon}
+c_2|\beta |^s+c_3,\quad\alpha ,\beta\in {\mathbb R},\label{e4}\end{equation}
holds. 
\par
Note that $F$ does not map $E$ 
into $E$. To address this difficulty we proceed as follows. 
Let $f_0(\xi )=f(\xi )-f(0)$ and let 
$F_0(x)(\xi )=f_0(x(\xi ))$. Then $f_0$ satisfies (\ref{e4}) and the 
mapping $F_0:E\to E$ is well 
defined. Equation (\ref{1}) can be rewritten in the form 
\[X(t)=S(t)x+\int_0^tS(t-s)F_0(X(s))ds+\int_0^tS(t-s)mds\]
\begin{equation}+\int_0^tS(t-s)\sqrt {Q}dW(s),\label{e41}\end{equation}
where $x\in E$ and $m(\xi )=f(0)$. Then 
\[\int_0^tS(t-s)mds\in\mbox{\rm dom}(A)\subset E,\]
and putting 
\[Z_m(t)=\int_0^tS(t-s)mds+\int_0^tS(t-s)\sqrt QdW(s),\]
and $Y(t)=X(t)-Z_m(t)$ we can rewrite (\ref{e41}) in the 
form 
\[Y(t)=S(t)x+\int_0^tS(t-s)F_0\left(Y(s)+Z_m(s)\right)ds.\]
Now it is clear, that the proof of existence and 
uniqueness of solutions provided in the proof of Theorem 
7.10 in \cite{dz1} applies in the present case. Moreover, 
\[\int_0^{\infty}\left\|S(t)m\right\|dt<\infty ,\]
and therefore the proof of Proposition \ref{t6}\ and,
consequently, all remaining statements in Section 2 remain 
valid as well. 
\par
Note that (\ref{e4}) is satisfied when $f$ is a 
polynomial of odd degree larger than one with a 
negative leading coefficient. It is well known (\cite{sin},
Theorem 2.2) that 
the subdifferential of the norm $\partial\left\|x\right\|$ at a point 
$x\in E$ contains the Dirac measures $\delta_{\xi_1}$ or $-\delta_{
\xi_2}$, if 
$\left\|x\right\|=x\left(\xi_1\right)$ or $\left\|x\right\|=-x\left
(\xi_2\right),$ respectively, hence it is 
easily seen that (\ref{e4}) implies (\ref{5}) and Hypothesis 1.3 is 
verified. The remaining assumptions depend on the 
covariance operator $Q$. Assume at first that $Q$ is 
boundedly invertible, that is, $Q^{}$ is an injection and 
$Q^{-1}\in\mathcal L(H).$ Then we have to verify (1.4) with $Q=I$ . By 
the well known estimates on the Green functions  
\cite{ar} it follows that 
$$||S(t)||_{HS}\le Ct^{-\frac{d}{4}},\quad t\in (0,1],$$
hence (1.4) is satisfied if the dimension $
d$ is 
one. In fact, for $d>1$ and $Q$ boundedly invertible even 
the Ornstein-Uhlenbeck process $Z$ does not take values in $H$, 
so these cases cannot be considered in the present 
framework. 
 If $d=1$, however, it is easy to see that 
all remaining conditions are satisfied. Proceeding as in Theorem 4.1 in 
\cite{peszat} we easily see that
$Z$ has an $E$-valued modification.
% and we have
%\begin{equation}\sup_{0\le\tau\le 1}\mathbb
%E\left\|\int_0^{\tau} S(\tau -r)Q^{1/2}d
%W(r)\right\|^p<\infty ,\label{e5}\end{equation}
%for each $p>0$. 
By the Sobolev embedding 
theorem, for each $\delta >\frac 14$ 
there exists a constant $c_{\delta}<\infty$ such that
\[||x||\le c_{\delta}|(-A)^{\delta}x|,\quad x\in\mbox{\rm dom}\left
(\left(-A_0\right)^{\delta}\right).\]
Take $\delta \in (\frac{1}{4} ,\frac{1}{2})$ and $p>2$; since
${\bf S}$ is exponentially stable
we obtain for some constants $c_1, c_2 ,c_3$ and $\omega >0$
\[\sup_{t\ge0}\mathbb E\left\|\int_0^tS(t-r)Q^{1/2}dW_r\right\|^p\]
%\[\le c_1\mathbb E\left\|\int_0^1S(r)Q^{1/2}dW_r\right\|^p+\sup_{
%t\ge 1}\mathbb E\left\|\int_1^tS(r)Q^{1/2}dW_r\right\|^p\]
\[\le c_1\sup_{t\ge 0}\mathbb E\left|\left(-A_0\right)^{\delta}
\int_1^tS(t-r)Q^{1/2}dW(r)\right|^p\le\]
\[\le c_2\sup_{t\ge 0}\left(\int_0^t\left\|\left(-A_0\right
)^{\delta}S(r)\right\|^2_{HS}dr\right)^{p/2}\le\]
%\[c_2+c_4\sup_{t\ge 0}\left(\int_0^{t-1}\left\|S(r)\right\|^2_{\mathcal 
%L(H)}\left\|\left(-A_0\right)^{\delta}S(1)\right\|^2_{HS}dr\right
%)^{p/2}\]
\[\le c_3\left(\int_0^{\infty}r^{-2\delta}e^{-2\omega r}dr\right)^{p/2}<\infty 
.\]
The same estimates hold for the process $Z_m$ and Hypothesis 1.2
is verified. Hypothesis 1.5 (strong Feller
property and topological irreducibility) is satisfied in the
present case (see e.g. \cite {MSe}) and we can conclude that in the
one-dimensional case if $Q$ is
boundedly invertible (in particular, if $Q=I$ which corresponds
to the case of space-time white noise) and the growth condition
(\ref{e4}) is satisfied Theorem \ref{t2} is applicable.
\par
Now we will examine some cases when the covariance $Q$ may
be degenerate. In order to obtain easily verifiable conditions we
only consider the so-called diagonal case. We assume that there 
exists an orthonormal basis $\left(e_n\right)$ in 
$H=L^2(D)$ such that $e_n\in E$ and for a certain $C<\infty$ 
 
\[\sup_{\xi\in D}|e_n(\xi )|<C,\quad\sup_{\xi\in D}|\nabla e_n(\xi 
)|<C\sqrt {\alpha_n},\quad n\ge 1,\]
and such that $(e_n),(\alpha_n)$ are the respective eigenvectors 
and eigenvalues of the operator $-A_0,$ $\alpha_n>\omega >0.$ We 
assume that the covariance operator $Q$ has the same 
eigenvectors $e_n$ with the respective eigenvalues 
$0<\lambda_n\le\lambda_0<\infty ,$ that is, 
\[Qe_n=\lambda_ne_n,\quad n\ge 1.\]
We again impose condition (\ref{e4}) on $f$. As in 
the previous case, we just have to check (\ref{2}) and Hypotheses 
1.2 and 1.5. For $\gamma\in (0,1)$ we have
\[\left\|S(t)Q^{1/2}\right\|^{2_{}}_{HS}=\sum_{n=1}^{\infty}\lambda_
ne^{-2\alpha_nt}\le\sup_{n\ge 1}\alpha_n^{1-\gamma}e^{-2\alpha_nt}
\sum_{n=1}^{\infty}\frac {\lambda_n}{\alpha_n^{1-\gamma}}\]
\[\le\frac {const}{t^{1-\gamma}}\sum_{n=1}^{\infty}\frac {\lambda_
n}{\alpha_n^{1-\gamma}}\]
thus (\ref{2}) is satisfied with $0<\delta <\gamma$ provided
\begin{equation}\sum_{n=1}^{\infty}\frac {\lambda_n}{\alpha_n^{1-
\gamma}}<\infty\label{e6}\end{equation}
holds for some $\gamma >0.$ Hypothesis 1.2 has been verified 
under condition (\ref{e6}) in Theorems 5.2.9 and 11.3.1 of 
\cite{dz2}. 
\par
The strong Feller property can be verified by Proposition
\ref{t10}. The assumption (\ref{52}) is obviously satisfied by
the local Lipschitz continuity of $f$. The condition (\ref{53})
is equivalent to 
 
\begin{equation}\sup_{n\in N}\left(\frac {\alpha_n}{\lambda_n}\left
(1-e^{-2\alpha_nt}\right)^{-1}\right)^{1/2}\in L^1\left(0,T \right
).\label{e9}\end{equation}
 It remains to verify that the solution to 
(\ref{e1}) is topologically irreducible. To this end, we will 
use Proposition 2.11 of \cite{M2} according to which, in the 
present case, it suffices to verify that $\mbox{\rm im}(\mathcal 
K)$ is dense 
in $C_0:=\{y\in C([0,T],E),y(0)=0\}$ where
\[\mathcal K:L^2(0,T,H)\to C_0,\quad\mathcal Ku(t):=\int_0^tS(t-r
)Q^{1/2}u(r)dr,\quad t\in [0,T].\]
The well known estimates on the Green 
kernel for a parabolic problem \cite{ar} yield
\begin{equation}\left\|S(t)\right\|_{\mathcal L(H,E)}\le\frac {co
nst}{t^{d/4}},\quad t\in [0,T],\label{e12}\end{equation}
hence $\|S(\cdot )\|_{\mathcal L(H,E)}$ is integrable for $d\le 3$ and the 
density follows from \cite{M2}, Proposition 2.8 and 
Remark 2.9. 
\par
We can summarise that Theorem 1.6 is applicable to 
the system (\ref{e1}) under conditions (\ref{e4}), (\ref{e6}), 
(\ref{e9}) provided 
$d\le 3$. In particular, if there exist $a\ge b\ge 0$ and 
constants $k_1,$$k_2$ such that
\[k_1\alpha_n^{-a}\le\lambda_n\le k_2\alpha_n^{-b},\quad n\ge 1,\]
or, equivalently, if
\[K_1n^{\frac {-2a}d}\le\lambda_n\le K_2n^{\frac {-2b}d},\quad n\ge 
1,\]
then it is easy to check that the condition (\ref{e6}) is 
satisfied if $b>\frac d2-1$ while (\ref{e9}) holds true if $a<1.$ So in 
this case our results are applicable if
\begin{equation}\frac d2-1<b\le a<1.\label{fin}\end{equation}

{\bf Acknowledgement.} The authors are grateful to Jan Seidler
for his valuable remarks.


\begin{thebibliography}{99}
\bibitem{aida}
Aida S.: Uniform positivity improving property, Sobolev 
inequalities and spectral gaps, {\em J. Funct. Analysis\/} 158 
(1998), 152-185
\bibitem{ar}
Arima R: On general boundary value problem for 
parabolic equations, {\em J.Math.Kyoto Univ.\/} 4(1964), 207-243
\bibitem{chen}
Chen Mu-Fa: Equivalence of exponential ergodicity and 
$L^2$-exponential convergence for Markov chains, 
{\em Stochastic Processes Appl.\/} 87 (2000), 281-297
\bibitem{pertex}
Chojnowska-Michalik A. and Goldys B.: Existence, 
uniqueness and invariant measures for stochastic semilinear 
equations in Hilbert spaces, {\em Probab. Theory Related }
{\em Fields\/} 102 (1995), 331-356
\bibitem{fock}
Chojnowska-Michalik A. and Goldys B.: 
Nonsymmetric Ornstein-Uhlenbeck semigroup as second 
quantized operator, {\em J. Math. Kyoto Univ.\/} 36 (1996), 481-498
\bibitem{dap1}
Da Prato G.: Large asymptotic behaviour of Kolmogorov equations in Hilbert 
spaces. Partial differential equations (Praha, 1998),
111-120, Chapman\& Hall/CRC Res. Notes Math., 406, Chapman\& Hall/CRC
\bibitem{dap2}
Da Prato Giuseppe: Poincar\'e  inequality for some measures in Hilbert 
spaces and application to spectral gap for transition semigroups, 
{\it Ann. Sci. Norm. Sup. Pisa} 25 (1997), 419--431
\bibitem{dz1}
Da Prato G. and Zabczyk J.: STOCHASTIC EQUATIONS IN 
INFINITE DIMENSIONS, Cambridge University Press 1992
\bibitem{dz2}
Da Prato G. and Zabczyk J.: ERGODICITY FOR INFINITE 
DIMENSIONAL SYSTEMS, Cambridge University Press 1996
\bibitem{dadego}
Da Prato G., Debussche A. and Goldys B.: 
Invariant measures of non symmetric dissipative
 stochastic systems, to appear in Probability Theory and Related Fields
\bibitem{dgz}
Da Prato G., G\c atarek D. and Zabczyk J.:  Invariant measures for 
semilinear stochastic equations, {\it Stochastic Anal. Appl.} 10
(1992), 387--408 
\bibitem{gg}
G\c atarek D. and Goldys B.: On invariant measures for 
diffusions on Banach spaces, {\em Potential Analysis\/} 7 (1997), 
539-553
\bibitem{gm}
Goldys B. and Maslowski B.: Ergodic control of semilinear 
stochastic equations and the Hamilton-Jacobi equation, {\em J. }
{\em Math. Analysis Appl.\/} 234 (1999), 592-631
\bibitem{JJ}
Jain N. and Jamison B: Contributions to Doeblin's theory 
of Markov processes, {\em Z.Wahrscheinlichkeitstheorie }
{\em Verw.Geb.\/} 8 (1967), 19-40
\bibitem{JR}
Jacquot S. and Royer G.: Ergodicite d'une classe 
d'equations aux derivees partielles stochastiques, 
{\em C.R.Acad.Sci.Paris Ser.Math.\/} 320 (1995), 231-236
\bibitem{La}
Lasota A. and Mackey M.C.: CHAOS, FRACTALS AND 
NOISE, Springer-Verlag, New York 1994
\bibitem{M1}
Maslowski B.: Strong Feller property for semilinear 
stochastic evolution equations and applications, 
STOCHASTIC SYSTEMS AND OPTIMISATION (Warsaw 1988), 
Lecture Notes in Control Inform. Sci 136, Springer 1989, 
210-224
\bibitem{M2}
Maslowski B.: On ergodic behaviour of solutions to 
systems of stochastic reaction-diffusion equations with 
correlated noise, STOCHASTIC PROCESSES and RELATED 
TOPICS (Georgenthal 1990), Akademie-Verlag, Berlin 1993, 
93-102
\bibitem{M3}
Maslowski B.: On probability distributions of solutions of 
semilinear stochastic evolution equations, {\em Stochastics }
{\em Stochastics Rep.\/} 45 (1993), 17-44
\bibitem{MSe}
Maslowski B. and Seidler J.: Probabilistic approach to 
the strong Feller property, {\em Probab. Theory Related }
{\em Fields\/} 118 (2000), 187-210
\bibitem{MSs}
Maslowski B. and Seidler J.: Invariant measures for 
nonlinear SPDE's: Uniqueness and stability, {\em Arch. Math. }
34 (1998), 153-172
\bibitem{meyn}
Meyn S. P. and and Tweddie R. L.: MARKOV CHAINS 
AND STOCHASTIC STABILITY, Springer-Verlag 1993
\bibitem{peszat}
Peszat S. and Seidler J.: Maximal inequalities and 
space-time regularity of stochastic convolutions, {\em Math. }
{\em Bohem.\/} 123 (1998), 7-32 
\bibitem{rr}
Roberts G. O. and Rosenthal J. S.: Geometric ergodicity 
and hybrid Markov chains, {\em Electron. Comm. Probab.\/} 2 
(1997), 13-25
\bibitem{roc}
R\"ockner M. and Zhang T. S.: Probabilistic representations and hyperbound 
estimates for semigroups, {\it Infin. Dimens. Anal. Quantum Probab.
Relat. Top.} 2 (1999), 337--358
\bibitem{Se1}
Seidler J.: Ergodic behaviour of stochastic parabolic 
equations, {\em Czechoslovak Math.J.\/} 47 (122) (1997), 
277-316
\bibitem{Sh}
Shardlow T.: Geometric ergodicity for stochastic PDEs, 
{\em Stochastic Anal.Appl.\/} 17 (1999), 857-869
\bibitem{sin}
Sinestrari E.: Accretive differential operators, {\em
Boll.Un.Mat.Ital B.\/} (5), 13(1976), 19-31 
\bibitem{wang}
Wang Feng-Yu: Functional inequalities, semigroup properties and spectrum 
estimates, {\it Infin. Dimens. Anal. Quantum Probab. Relat. Top.} 3
(2000), 263--295.
\bibitem{wu}
Wu Liming: Uniformly integrable operators and large 
deviations for Markov processes, {\em J. Funct. Analysis\/} 172 
(2000), 301-376
\end{thebibliography}
\end{document}